\documentclass{article}
\usepackage{amsmath}
\usepackage{amsthm}
\usepackage{amsfonts}
\usepackage{amssymb}
\usepackage{tikz}
\usepackage{hyperref}
\usepackage[margin=1in]{geometry}

\theoremstyle{theorem}
\newtheorem{theorem}{Theorem}
\newtheorem{proposition}[theorem]{Proposition}
\newtheorem{lemma}[theorem]{Lemma}
\newtheorem{corollary}[theorem]{Corollary}

\theoremstyle{definition}

\DeclareMathOperator{\tr}{\texttt{R}}
\DeclareMathOperator{\tl}{\texttt{L}}
\DeclareMathOperator{\C}{\mathcal{C}}
\DeclareMathOperator{\R}{\mathcal{R}}
\DeclareMathOperator{\Lf}{\mathcal{L}}
\DeclareMathOperator{\I}{\mathcal{I}}
\DeclareMathOperator{\Ou}{\mathcal{O}}
\DeclareMathOperator{\pr}{pr}
\DeclareMathOperator{\E}{\mathbb{E}}
\DeclareMathOperator{\V}{\mathbb{V}}
\DeclareMathOperator{\nap}{\curlywedge}

\begin{document}

\title{The Most Malicious Maitre d'}

\markright{Abbreviated Article Title}
\author{Tejo Madhavarapu and T. Kyle Petersen}  

\maketitle

\begin{abstract}
We revisit the family of ``napkin problems'' first discussed by Peter Winkler in his 2004 puzzle book and discussed subsequently in several papers. This family of problems involves diners being seated around a circular table in which napkins are placed between the seats, so it is ambiguous which napkin belongs to which seat. The problems seek to quantify the proportion of diners who end up with no napkin under a variety of assumptions, one of which is that the diners have an adaptive adversary known as ``the malicious maitre d','' who seeks to maximize the number of napkinless diners. Winkler's original strategy for the adaptive maitre d' was shown to be suboptimal in 2023 when Acton, Petersen, Shirman, and Toal presented a better (yet also suboptimal) strategy. In this paper we demonstrate an optimal strategy for the malicious maitre d' and compute the expected proportion of napkinless diners as a function of the probability of a diner taking the left napkin. Interestingly, the strategy does not depend on the exact probability; rather, it only depends on which napkin (left or right) the diners tend to prefer.
\end{abstract}


\section{Napkin problems}

Peter Winkler's book \textit{Mathematical Puzzles: A Connoisseur's Collection} \cite{Winkler} contained a pair of ``napkin problems'' which have inspired several papers over the past couple of decades \cite{Sudbury, CP, Eriksen, APSToal, APSTenner}. In both problems there is  a large circular table with napkins placed between the seats, so that it is unclear which napkin belongs with which seat. There are $n$ seats, $n$ napkins, and we have $n$ diners about to sit at the table. Each diner sits at a distinct time step and attempts to take a napkin. If both the left and right napkins are available, they choose randomly. If only one napkin is available, they take that napkin. If neither napkin is available the diner is \textbf{napkinless}.

The first version of the problem, posed by John Conway, asks for the expected proportion of napkinless diners when the diners sit randomly. The second version of the problem, that Winkler attributes to Rob Pike, asks for a worst-case analysis. This problem imagines an adversary, \textbf{the Malicious Maitre d'}, who tells the diners where to sit, and whose goal is to maximize the expected proportion of napkinless diners. The maitre d' can operate under different assumptions. In the first version from Winkler's book \cite[p. 22]{Winkler}, the maitre d' is devising a malicious seating order. The diners arrive one by one, and this maitre d' tells each diner where to sit. But also discussed by Winkler is \textbf{the Adaptive Maitre d'} \cite[p. 25]{Winkler}. Here Winkler introduces a version of the question in which the maitre d' can witness the napkin selected by each diner as they sit down and use that information to decide where to place subsequent diners. A final variation studied in \cite{APSTenner} imagines a \textbf{Clairvoyant Maitre d'} who learns what napkin a diner will choose just before they sit down. The reader is encouraged to experiment with these variants on \cite{website}, a web application created for exploring this problem.

Here is a quick rundown of known results so far. For large tables, the expected proportion of napkinless diners is:
\begin{itemize}
\item $(2-\sqrt{e})^2 \approx 12.34\%$ in Conway's napkin problem,
\item $9/64 \approx 14.06\%$ for Pike's problem,
\item $1/6 \approx 16.67\%$ for Winkler's solution to the Adaptive Maitre d' problem,
\item between $17.72\%$ and $18.15\%$ using the solution to the Adaptive Maitre d' problem in \cite{APSToal}, and
\item $1/3 \approx 33.33\%$ for the solution to the Clairvoyant Maitre d' problem, which was proved to be optimal in \cite{APSTenner}.
\end{itemize}
All of the above results are stated under the assumption that diners have an equal probability of choosing the left napkin or the right napkin. In the case of Conway's problem, however, each of the papers \cite{Sudbury, CP, Eriksen} generalize to the case of a fixed probability $p$ of taking the left napkin, probability $q=1-p$ of taking the right napkin. The expectation $(2-\sqrt{e})^2$ result is generalized to
\begin{equation}\label{eq:conway}
 \frac{(1-pe^q)(1-qe^p)}{pq}.
\end{equation}
See \cite[Theorem 6]{Sudbury}, \cite[Theorem 4]{CP}, and \cite[Section 1]{Eriksen}.

The contribution of this paper is to settle the question of an optimal strategy for the Adaptive Maitre d'. Winkler initially claimed his strategy for this problem, dubbed ``trap setting'' in \cite{APSToal}, was optimal. Here is that strategy in Winkler's words \cite[p.25]{Winkler}:
\begin{quote}
 If the first diner takes (say) his right napkin, the next is seated two spaces to his right so that the diner in between may be trapped. If the second diner also takes his right napkin, the maitre d' tries again by skipping another chair to the right. If the second diner takes his left napkin (leaving the space between him and the first diner napkinless), the third diner is seated directly to the second diner's right. Further diners are seated according to the same rule until the circle is closed, then the remaining diners (some of whom are doomed to be napkinless) are seated.
\end{quote}

The thrust of \cite{APSToal} was to carefully analyze trap setting and to prove their strategy of ``napkin shunning'' was superior. Here is that strategy according to \cite[p.729]{APSToal}:
\begin{quote}
 The insight for the new strategy begins with the observation that there is a bijection between napkins at the table and diners at the table. Hence for each napkinless diner, there must be, somewhere on the table, a dinerless napkin ignored by both of its adjacent diners. Rather than directly attempting to set traps by having diners sit apart and hoping they reach toward each other, the napkin shunning strategy seats a new diner on the side away from the napkin the previous diner has chosen, in the hopes that the second diner reaches away from the first diner, shunning the napkin between. Once we succeed in shunning a napkin, the next diner is seated in the middle of a gap of empty seats and we attempt to shun another napkin. We continue to split gaps and try to shun napkins until all gaps are filled. Since every shunned napkin will ultimately correspond to a napkinless diner, the traps will have set themselves.
\end{quote}

\begin{figure}
\includegraphics[width=15cm]{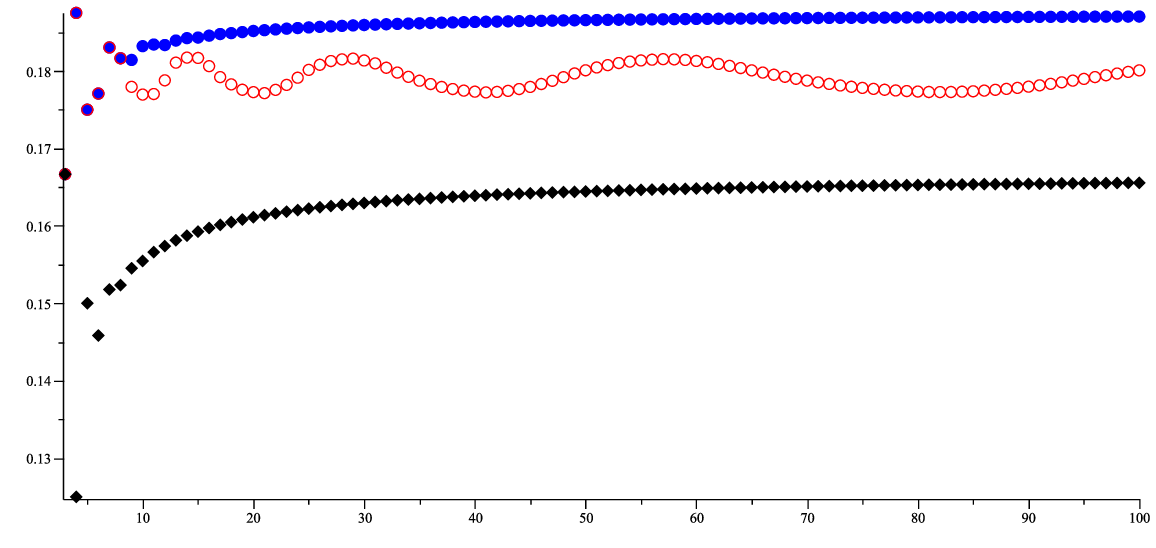}
\caption{The expected proportion of napkinless diners for $p = 1/2$, when using trap-setting (black diamonds), napkin shunning (red open circles), and napkin shunning with wide traps (blue solid circles).}\label{fig:comps}
\end{figure}

The authors of \cite{APSToal} knew that while their strategy was better than trap setting, it was not optimal. They asked if an optimal strategy exists, and this paper proves the answer is yes, via a strategy that optimizes the choice of where to place a diner once a napkin has been successfully shunned. Moreover, this strategy is optimal for any fixed probability $p$, not only the $p=q=1/2$ case. It turns out that the ``split the gap'' approach described above agrees with the optimal strategy for tables with at most eight seats, but in reality the best choice is to place the next diner three seats away from one of the diners who shunned the napkin. This is similar to trap setting, but with two empty seats instead of one. We find this solution, which we call \textbf{napkin shunning with wide traps} to be beautifully mundane. We compare the performance of the three algorithms (trap setting, napkin shunning, and napkin shunning with wide traps) for $p=q=1/2$ and $3\leq n\leq 100$ in Figure \ref{fig:comps}.

\begin{theorem}\label{thm:main}
For the adaptive maitre d' problem, napkin shunning with wide traps maximizes the expected number of napkinless diners. As the table size gets large, the expected proportion of napkinless diners converges to 
\begin{equation}\label{eq:expr}
 pq\left(\frac{1+r}{1+2r}\right) = \frac{r(1-r^2)}{1+2r},
\end{equation}
where $r = \min\{p,q\}$.
\end{theorem}

At $p=q = 1/2$, this expectation is $3/16=18.75\%$. However, unlike Equation \eqref{eq:conway} for the case for Conway's napkin problem, the expectation in \eqref{eq:expr} is \emph{not} maximized at this value of $p$. Rather, its maximum occurs at $r \approx 0.4554$ (the real root of $1-3r^2-4r^3$), when the expectation is that roughly $18.89\%$ of diners are napkinless.  So the most malicious maitre d' actually hopes for mild agreement amongst the guests! See the solid red curve in Figure \ref{fig:rgraph}. 

\begin{figure}
\[
\begin{tikzpicture}[xscale=10,yscale=30]
\draw (0,-0.01)--(0,0.2);
\draw (-0.05,0)--(1.05,0) node[below right] {$p$};
\foreach \x in {1,...,5}{
 \draw (-0.005,-0.02+0.04*\x)--(0,-0.02+0.04*\x);
 \draw (-0.1+0.2*\x,-0.0025)--(-0.1+0.2*\x,0);
}
\foreach \x in {0.04,0.08,0.12,0.16,0.20}{
\draw (-0.01,\x) node[left] {\footnotesize{$\x$}} --(0,\x);
}
\foreach \x in {0.2,0.4,0.6,0.8,1}{
\draw (\x,-0.005) node[below] {\footnotesize{$\x$}} --(\x,0);
}
\draw[dotted] (0.5,-0.01)--(0.5,0.2); 
\draw[green!80!black, domain=0.001:0.5] plot[smooth] (\x,{\x*(1-\x)/(1+\x)});
\draw[green!80!black, domain=0.0001:0.5] plot[smooth] (1-\x,{\x*(1-\x)/(1+\x)});
\draw[red, domain=0.001:0.5] plot[smooth] (\x,{\x*(1-\x^3)/(1+3*\x)});
\draw[red, domain=0.001:0.5] plot[smooth] (1-\x,{\x*(1-\x^3)/(1+3*\x)});
\draw[blue, domain=0.001:0.999] plot[smooth] (\x,{(1-\x*exp(1-\x))*(1-(1-\x)*exp(\x))/(\x*(1-\x))});
\draw[very thick, domain=0.0001:0.5] plot[smooth] (\x,{\x*(1-\x^2)/(1+2*\x)});
\draw[very thick, domain=0.0001:0.5] plot[smooth] (1-\x,{\x*(1-\x^2)/(1+2*\x)});
\draw[thin] (0.29,0.16)--(0.1,0.17) node[draw=black,fill=white,inner sep=2] {\footnotesize{Width $1$ traps}};
\draw[thin] (0.755,0.14)--(0.92,0.15) node[draw=black,fill=white,inner sep=2] {\footnotesize{Width $3$ traps}};
\draw[thin] (0.4554,0.1889)--(0.5,0.21) node[draw=black,fill=white,inner sep=2] {\footnotesize{Width $2$ traps (optimal)}};
\draw[thin] (0.5,0.1234)--(.4,0.14) node[draw=black,fill=white,inner sep=2] {\footnotesize{Random}};
\end{tikzpicture}
\]
\caption{The expected proportion of napkinless diners as a function of $p$. The thick black curve is the optimal expectation for the adaptive ma\^itre d', coming from napkin shunning with traps of width 2. The green curve corresponds to traps of width 1 and the red curve corresponds to traps of width 3. The expectation for Conway's napkin problem, which is equivalent to the ma\^itre d' following a random strategy, is shown in blue. 
}
\label{fig:rgraph}
\end{figure}

Beyond settling the question of the adaptive maitre d', we believe this work is also an advertisement for the necessity of thorough computational experimentation. As mentioned, napkin shunning with wide traps agrees with the napkin shunning strategy of \cite{APSToal} on tables of up to $n=8$ seats, but a greedy search for optimal seating strategies makes the best strategy clear after investigating the next few cases.

We are reminded of the old adage from the early days of computer science: ``premature optimization is the root of all evil'' and suggest that for this problem ``premature generalization is the root of all evil.''

\section{Precise statement of problem and solution}

We follow \cite{APSToal} for much of our notation and terminology. We fix a positive integer $n$ to represent the number of diners, and we fix $p$ for the probability that a diner will reach to the left and $q=1-p$ for the probability that they reach right. We assume this probability is identical and independent for each diner. For the bulk of the paper, we choose to assume $p\leq 1/2$ to make our arguments simpler. This is no loss of generality, since if $p>1/2$, all the arguments are mirrored by swapping ``left'' with ``right'' and $p$ with $q$. We will make clear when we wish to consider $p>1/2$.

A key data structure for us is the sequence of diner preferences for their right napkin or left napkin, which we call a \textbf{preference order}. A preference order for $n$ diners is a list $\sigma \in \{\tl,\tr\}^n$. That is, $\sigma = (\sigma_1, \sigma_2, \ldots, \sigma_n)$, where $\sigma_j = \tl$ means the $j$th diner prefers the napkin to their left and $\sigma_j = \tr$ means the diner prefers the napkin to their right. A preference order $\sigma$ occurs with probability $p^kq^{n-k}$, where $k$ is the number of ``$\tl$'' entries of $\sigma$.

A \textbf{strategy} is a rule that assigns a seat number to a diner as they arrive. By the symmmetry of the circular table, the first diner's assignment is unimportant, and we distinguish strategies only after the first diner has been seated. From that point forward, the input for the strategy is the current state of the table, i.e., which seats are occupied and which napkins are still available. 

If there are $n$ seats at the table, we label positions from left to right from the perspective of the diners (counterclockwise if viewing the table from above) and the diners are labeled by their order of arrival at the table. We choose to declare seat ``1'' as the seat in which the first diner sits. Thus as we move to the right from diner 1, we encounter seats 2, 3, and so on. In moving to the left of diner 1 we find seats $n$, $n-1$, and so on. We write a \textbf{seating order} as a list $w=(w_1,w_2,\ldots,w_n)$, where $w_i = j$ means that the $j$th diner sits in seat $i$. In principle, the maitre d' could create a seating order for any permutation $w$ of $\{1,2,\ldots,n\}$ such that $w_1=1$.

When providing examples, we will often draw tables as arrays with the seating order $w$ on bottom and the corresponding diner preference above ($\tr$ or $\tl$). To emphasize when a diner does not receive their desired napkin, we write $\tl'$ or $\tr'$ for diners who get a napkin other than their preferred napkin and we write $\tl^*$ or $\tr^*$ for napkinless diners. In Figure \ref{fig:circle} we see two visualizations of the seating order $w=(1,5,2,8,4,6,7,3)$ together with preference order $\sigma=(\tr,\tl,\tl,\tr,\tr,\tl,\tr,\tl)$. Notice that the sequence of events here conspire to leave two of the diners (diners 5 and 7) without a napkin.

\begin{figure}
\[
\begin{array}{ccc}
\begin{tikzpicture}[scale=1, baseline=0cm]
\draw (0,0) circle (2);
\draw[dashed] (360/8*5+90:2.25)--(360/8*4+112.5:1.75);
\draw[thick] (360/8*0+90:2.25) node[fill=white,inner sep=0] {$1$} -- (360/8*0+112.5:1.75) node[circle,fill=white,inner sep=1] {$\curlywedge$};
\draw (360/8*1+90:2.25) node[fill=white,inner sep=1, draw=black, circle] {$5$};
\draw[thick] (360/8*2+90:2.25) node[fill=white,inner sep=0] {$2$} -- (360/8*1+112.5:1.75) node[circle,fill=white,inner sep=1] {$\curlywedge$};
\draw[thick] (360/8*3+90:2.25) node[fill=white,inner sep=0] {$8$} -- (360/8*2+112.5:1.75) node[circle,fill=white,inner sep=1] {$\curlywedge$};
\draw[thick] (360/8*4+90:2.25) node[fill=white,inner sep=0] {$4$} -- (360/8*4+112.5:1.75) node[circle,fill=white,inner sep=1] {$\curlywedge$};
\draw[thick] (360/8*5+90:2.25) node[fill=white,inner sep=0] {$6$} -- (360/8*5+112.5:1.75) node[circle,fill=white,inner sep=1] {$\curlywedge$};
\draw (360/8*6+90:2.25) node[fill=white,inner sep=1,draw=black,circle] {$7$};
\draw[thick] (360/8*7+90:2.25) node[fill=white,inner sep=0] {$3$} -- (360/8*6+112.5:1.75) node[circle,fill=white,inner sep=1] {$\curlywedge$};
\draw (360/8*3+112.5:1.75) node {$\curlywedge$};
\draw (360/8*7+112.5:1.75) node {$\curlywedge$};
\end{tikzpicture}
&
\leftrightarrow
&
 \begin{array}{cccccccc}
   \tr&\tr^*&\tl&\tl&\tr&\tl'&\tr^*&\tl \\
   \hline
   1&{\color{red} \mathbf{5}}&2&8&4&6&{\color{red} \mathbf{7}}&3\\
 \end{array}
\end{array}
\]
\caption{A seating arrangement for eight, with preferences that lead to two napkinless diners. Notice that diner $6$ preferred the left napkin but took the right napkin because the left was not available.}\label{fig:circle}
\end{figure}

We now define our optimal strategy via a deterministic algorithm that we call \textbf{Algorithm B} (``B'' for best). Here the maitre d' makes a choice about where to seat diner $i$ based on the observed choices of diners $1, \ldots, i-1$.

\subsection{Algorithm B (napkin shunning with wide traps).}

We initialize by seating the first diner and observing which napkin they select. The algorithm then proceeds as follows.
\begin{enumerate}
\item[\textbf{B1.}] (Napkin shunning) If all diners are seated, stop. If not, find the leftmost napkin (relative to seat 1) that neighbors only one empty seat. Place a diner in that seat and repeat step B1. If there is no such napkin, go to step B2.
\item[\textbf{B2.}] (Wide trap setting) Place the next diner 2 seats from the right edge of the leftmost gap, setting a wide trap, if this gap has at least three empty seats. If the gap has fewer than three empty seats, place the new diner in the leftmost empty seat of this gap. Return to step B1. 
\end{enumerate}

We remark that there is some flexibility in the definition of step B2, e.g., the choice of gap does not matter for the expected performance of the algorithm. We specify a particular gap for concreteness.

We write $\nu(\sigma)$ to denote the number of napkinless diners that result from applying Algorithm B to preference order $\sigma$. For example, suppose the preferences of 18 diners are given by the preference order below:
\[
 \sigma = (\tr,\tr,\tl,\tr,\tl,\tl,\tr,\tl,\tl,\tr,\tr,\tl,\tl,\tr,\tr,\tl,\tl,\tr).
\]

After the first twelve diners have been seated, we see
\[
\begin{array}{cccccccccccccccccc}
   \tr&-&\tl&\tr&-&-&\tl&\tl&\tr'&\tl&\tr&-&\tl&\tr&-&-&\tl&\tr \\
   \hline
   1&&12&11&&&8&9&10&6&7&&5&4&&&3&2\\
 \end{array}
\]
and now we simply fill in from left to right to get our final seating arrangement. None of these first twelve diners are napkinless (though diner 10 is frustrated). However, we shunned four napkins, and therefore the last person to sit in each of the remaining gaps will be napkinless. After all the diners are seated, we see:
\[
\begin{array}{cccccccccccccccccc}
   \tr&\tl^*&\tl&\tr&\tr&\tr^*&\tl&\tl&\tr'&\tl&\tr&\tl^*&\tl&\tr&\tl'&\tr^*&\tl&\tr \\
   \hline
   1&{\color{red} \mathbf{13}}&12&11&14&{\color{red} \mathbf{15}}&8&9&10&6&7&{\color{red} \mathbf{16}}&5&4&17&{\color{red} \mathbf{18}}&3&2\\
 \end{array}
\]

\subsection{Table types}

When searching for an optimal strategy, one must compare all choices that could be made for where to place the next diner at the table. The available seats come in intervals separated by seated diners. We refer to these runs of open seats as ``table types'' and classify them as follows. Each table type is indexed by the total number of empty seats, $n$.

\begin{itemize}

\item \textbf{The right-leaning table.} This table, denoted $\R_n$, is an interval with $n$ empty seats and $n$ napkins, such that the leftmost seat has no napkin to its left, e.g.,
\[
\begin{tikzpicture}
\draw (0,0) node[left] {$\R_7=$};
\draw (0,0) node[right, draw=black] {
\begin{tikzpicture}
    \foreach \x in {1,...,7}{
    \draw (\x-.25,0)--(\x+.25,0);
    \draw (\x+.2,.5) node {$\nap$};
    }
\end{tikzpicture}
};
\end{tikzpicture}
\]
It arises from an interval of the form $\tr-\cdots -\tr$ or $\tl'-\cdots -\tr$.

\item \textbf{The left-leaning table.} This table, denoted $\Lf_n$, is an interval with $n$ empty seats and $n$ napkins, such that the rightmost seat has no napkin to its right, e.g.,
\[
\begin{tikzpicture}
\draw (0,0) node[left] {$\Lf_7=$};
\draw (0,0) node[right, draw=black] {
\begin{tikzpicture}
    \foreach \x in {1,...,7}{
    \draw (\x-.25,0)--(\x+.25,0);
    \draw (\x-.8,.5) node {$\nap$};
    }
\end{tikzpicture}
};
\end{tikzpicture}
\]
It corresponds to an interval of the form $\tl-\cdots -\tl$ or $\tl-\cdots -\tr'$.

\item \textbf{The inner-facing table.} This table, denoted $\I_n$, is an interval with $n$ empty seats and $n-1$ napkins, such that the leftmost seat has no napkin to its left and the rightmost seat has no napkin to its right, e.g.,
\[
\begin{tikzpicture}
\draw (0,0) node[left] {$\I_7=$};
\draw (0,0) node[right, draw=black] {
\begin{tikzpicture}
    \foreach \x in {1,...,6}{
    \draw (\x-.25,0)--(\x+.25,0);
    \draw (\x-.8,.5) node {$\nap$};
    }
    \draw (0.25,0)--(-0.25,0);
\end{tikzpicture}
};
\end{tikzpicture}
\]
It corresponds to an interval of the form $\tr-\cdots -\tl$, $\tl'-\cdots -\tl$, $\tr-\cdots-\tr'$, or $\tl'-\cdots - \tr'$.

\item \textbf{The outer-facing table.} This table, denoted $\Ou_n$, is an interval with $n$ empty seats and $n+1$ napkins, such that each seat has a napkin to both its left and its right, e.g.,
\[
\begin{tikzpicture}
\draw (0,0) node[left] {$\Ou_7=$};
\draw (0,0) node[right, draw=black] {
\begin{tikzpicture}
    \foreach \x in {1,...,7}{
    \draw (\x-.25,0)--(\x+.25,0);
    \draw (\x-.8,.5) node {$\nap$};
    }
    \draw (7.2,.5) node {$\nap$};
\end{tikzpicture}
};
\end{tikzpicture}
\]
It corresponds to an interval of the form $\tl-\cdots -\tr$.

\end{itemize}

Within a given table type, we label the seats $1,2,\ldots, n$ from left to right.

For completeness, we also let $\C_n$ denote a circular table with $n$ empty seats and $n$ napkins between the seats, prior to any diners being seated. We also remark that a careful reader might notice that Algorithm B will never produce an interval of type $\Ou_k$ when starting from $\C_n$. However, to prove optimality of Algorithm B, we need to compare it to strategies that might produce such an interval.

\subsection{Discovering an optimal strategy}

The first diner to sit at the table creates a table of type $\Lf_{n-1}$ or of type $\R_{n-1}$, depending on the napkin they choose. Any strategy that that is optimal for $\C_n$ must therefore be able to identify the optimal seat on one of these smaller linear tables. What happens when a next seat is chosen? 

Let $\E(\mathcal{T}_n)$ denote the maximal expectation for the number of napkinless diners on a table of type $\mathcal{T}_n$, considered over all possible strategies. If there is only one seat, there is only one seating strategy (put them in the seat). Thus we observe $\E(\I_1) = 1$ (since there is no napkin) and $\E(\Lf_1) = \E(\R_1) = \E(\Ou_1) = 0$. We also define $\E(\mathcal{T}_0) = 0$ for convenience.

Now consider an interval of type $\Lf_k$ with $k>1$. If we place the next diner in seat 1 of this interval, we end up with an interval of type $\Lf_{k-1}$ if the diner reaches left (with probability $p$) or an interval of type $\I_{k-1}$ if the diner reaches right (with probability $q$). If we place the diner in seat $k$, that diner has no choice but to take their left napkin and we end up with an interval of type $\Lf_{k-1}$. If we place the diner in some seat $1<i<k$, we break the interval into two disjoint intervals. If the diner takes their left napkin we have created an interval of type $\Lf_{i-1}$ and one of type $\Lf_{k-i}$. If the diner takes the right napkin, we break the interval into an interval of type $\Ou_{i-1}$ and of type $\I_{n-i}$. We illustrate the idea with $n=9$ and $i=4$ in Figure \ref{fig:Lsplit}.

\begin{figure}
\[
\begin{tikzpicture}
\draw (0,0) node[scale=.8]{
\begin{tikzpicture}[>=stealth]
\draw (0,0) node[left] {$\Lf_9=$};
\draw (0,0) node[right, draw=black, scale=.7] (a) {
\begin{tikzpicture}
    \foreach \x in {1,...,9}{
    \draw (\x-.25,0)--(\x+.25,0);
    \draw (\x-.8,.5) node {$\nap$};
    }
    \draw (3.8,.25) node {?};
\end{tikzpicture}
};
\draw (-4,-2) node[right,draw=black,scale=.7] (b) {
\begin{tikzpicture}
    \foreach \x in {1,...,9}{
    \draw (\x-.25,0)--(\x+.25,0);
    \draw (\x-.8,.5) node {$\nap$};
    }
    \draw (3.75,0) node[fill=white] {$\tl$};
    \draw[line width=10,gray,opacity=.3,cap=round] (4,0)--(3.5,.5);
\end{tikzpicture}
};
\draw (4,-2) node[right,draw=black,scale=.7] (c) {
\begin{tikzpicture}
    \foreach \x in {1,...,9}{
    \draw (\x-.25,0)--(\x+.25,0);
    \draw (\x-.8,.5) node {$\nap$};
    }
    \draw (3.75,0) node[fill=white] {$\tr$};
    \draw[line width=10,gray,opacity=.3,cap=round] (4,0)--(4.5,.5);
\end{tikzpicture}
};
\draw (-2,-3) node[right] {$\Lf_3 \cup \Lf_5$};
\draw (6,-3) node[right] {$\Ou_3 \cup \I_5$};
\draw[->] (a)-- node[midway,left,fill=white,inner sep=1] {$p$} (b);
\draw[->] (a)-- node[midway,right,fill=white,inner sep=1] {$q$} (c);
\end{tikzpicture}
};
\end{tikzpicture}
\]
\caption{Placing a diner into the middle of one table type produces a union of two types of tables of smaller sizes.}\label{fig:Lsplit}
\end{figure}

Thus to choose a seat that maximizes the expected number of napkinless diners for a fixed $k$, we can compute recursively over these cases. That is, the maximum expectation for $\Lf_k$ is
\begin{equation} \label{eq:ELmax} \E(\Lf_k) = \max\{\E(\Lf_{k-1}), \; p(\E(\Lf_{i-1}) + \E(\Lf_{k-i})) + q(\E(\Ou_{i-1}) + \E(\I_{k-i})) \; :\; 1 \leq i < k \}.
\end{equation} 
To compute this maximum, we need a way to find the maximum expectation on smaller intervals of types $\Lf$, $\Ou$, and $\I$. But from investigating $\Ou$ and $\I$ we are inevitably led to consider intervals of type $\R$ as well. By considering each in turn, we have:
\begin{align}
\E(\R_k) &= \max\{
    \E(\R_{k-1}),\; 
    p(\E(\I_{i-1}) + \E(\Ou_{k-i})) 
    + q(\E(\R_{i-1}) + \E(\R_{k-i})) 
    : 1 < i \leq k
\}, \label{eq:ERmax} \\
\E(\Ou_k) &= \max\{
    p(\E(\Lf_{i-1}) + \E(\Ou_{k-i})) 
    + q(\E(\Ou_{i-1}) + \E(\R_{k-i})) 
    : 1 \leq i \leq k
\}, \label{eq:EOmax} \\
\E(\I_k) &= \max\{
    \E(\I_{k-1}),\; 
    p(\E(\I_{i-1}) + \E(\Lf_{k-i})) 
    + q(\E(\R_{i-1}) + \E(\I_{k-i})) 
    : 1 < i < k
\}. \label{eq:EImax}
\end{align}
For the sake of completeness, we also record our initial observation of this subsection, i.e.,
\begin{equation} \label{eq:ECmax}
\E(\C_{k}) = p\E(\Lf_{k-1}) + q\E(\R_{k-1}).
\end{equation} 

Using these recurrences along with our base cases of $k=0$ and $k=1$ allows us to compute the first few maximal expectations. These maximal expectations are listed in Table \ref{tab:smallcases}.

\begin{table}
\[
\begin{array}{r| c | c | c | c | c}
 n & \E(\I_n) & \E(\Lf_n) & \E(\R_n) & \E(\Ou_n) & \E(\C_n)\\
 \hline
 1 & 1 & 0 & 0 & 0 & 0 \\
 2 & 1 & q & p & 0 & 0 \\
 3 & 1 & q + pq & p + pq & pq & 2pq \\
 4 & 1 + pq & q + pq(1 + p) & p + pq(1 + q) & 2pq & 3pq \\
 5 & 1 + 2pq & q + pq(2 + p^2) & p + pq(2 + q^2)  & pq(3-pq) & pq(4-2pq)
\end{array}
\]
\caption{Maximal expectations for small table types.}\label{tab:smallcases}
\end{table}

Along with our recursive computations of maximal expectation, we can also have our computer tell us which seat placements achieve this maximum. On tables of type $\Lf_k$, $\R_k$, and $\Ou_k$, we find confirmation for napkin shunning, as the best placement occurs on an end of the table where there is more than one napkin available, i.e., seat 1 on $\Lf_k$, seat $k$ on $\R_k$, and either end seat for table $\Ou_k$. However, for table $\I_k$ we start to see something different than in previously studied strategies.

With $k\leq 6$, the maximum for $\E(\I_k)$ is achieved by placing the next diner in the middle seat if $k$ odd, and just to the right or left of middle when $k$ is even (depending on $p$). This agrees with the strategy for napkin shunning outlined in \cite{APSToal}, and is verifiably superior to trap setting. But starting with $k=7$ we get something different! 

For $\E(\I_7)$ consider the maximum of $p(\E(\I_{i-1}) + \E(\Lf_{7-i})) + q(\E(\R_{i-1}) + \E(\I_{7-i}))$ over $i=1,2,\ldots, 7$. This is highlighted in bold in Table \ref{tab:I7} for several values of $p\leq 1/2$. We see that the maximum is indeed found in seat $k-2 = 5$. (We find the maximum in seat $3$ if $p\geq 1/2$.) In contrast, the strategy of \cite{APSToal} would recommend placing the next diner in seat $4$, which is only the fifth best choice!

\begin{table}
\[
\begin{array}{r|lllll}
 i & p=0.1 & p=0.2 & p=0.3 & p=0.4 & p=0.5 \\
 \hline
 1 & 1.261 & 1.448 & 1.567 & 1.624 & 1.625  \\
 2 & 1.27009 & 1.48128 & 1.63567 & 1.73536 & 1.78125  \\
 3 & 1.2709 & 1.4864 & 1.6489 & 1.7584 & \mathbf{1.8125}  \\
 4 & 1.27 & 1.48 & 1.63 & 1.720 & 1.75  \\
 5 & \mathbf{1.3429} & \mathbf{1.5824} & \mathbf{1.7329} & \mathbf{1.8064} & \mathbf{1.8125}  \\
 6 & 1.33561 & 1.56192 & 1.70203 & 1.77184 & 1.78125  \\
 7 & 1.261 & 1.448 & 1.567 & 1.78125 & 1.625 
\end{array}
\]
\caption{Identifying optimal seats for $\I_7$ and various values of $p$.}\label{tab:I7}
\end{table}

\section{Recurrences, generating functions, and limiting expectation for Algorithm B}

Before we prove Algorithm B is optimal, it will be helpful to develop some enumerative properties first. Write $\E_B(\mathcal{T}_n)$ for the expected number of napkinless diners on table $\mathcal{T}_n$ using Algorithm B. Our ultimate goal will be to show that $\E_B(\mathcal{T}_n) = \E(\mathcal{T}_n)$, and to do this we develop useful formulas for $\E_B(\mathcal{T}_n)$.

The application of Algorithm B to each table type is straightforward. By the linearity of expectation we achieve the following recurrences.

\begin{proposition}[Recurrences for expectations]\label{prp:Erecs}
We have $\E_B(\C_1) = \E_B(\Ou_1) = \E_B(\Lf_1) = \E_B(\R_1) = 0$, and for $n\geq 2$:
 \begin{align}
 \E_B(\C_n) &= p\E_B(\Lf_{n-1})+q\E_B(\R_{n-1}), \label{eq:Crec} \\
 \E_B(\Ou_n) &= p\E_B(\Ou_{n-1}) + q\E_B(\R_{n-1}), \label{eq:Orec} \\
 \E_B(\Lf_n) &= p\E_B(\Lf_{n-1}) + q\E_B(\I_{n-1}), \label{eq:Lrec} \\
 \E_B(\R_n) &= p\E_B(\I_{n-1}) + q\E_B(\R_{n-1}). \label{eq:Rrec}
 \end{align}
Further, $\E_B(\I_1) = \E_B(\I_2) = \E_B(\I_3) = 1$, and for $n\geq 4$: 
  \begin{equation}\label{eq:Irec}
   \E_B(\I_n) =  \E_B(\R_{n-2}) + \E_B(\Lf_3) .
  \end{equation}
\end{proposition}

\begin{proof}
Equation \eqref{eq:Crec} follows no matter what seat the first diner takes at the circular table. 

On a table of type $\Lf_n$, the next diner is placed in seat 1, resulting in a table of type $\Lf_{n-1}$ with probability $p$ and a table of type $\I_{n-1}$ with probability $q$. On a table of type $\Ou_n$, the next diner is placed in seat 1, resulting in a table of type $\Ou_{n-1}$ with probability $p$ and a table of type $\R_{n-1}$ with probability $q$. On a table of type $\R_n$, the diner is placed in seat $n$, which gives $\I_{n-1}$ with probability $p$ and $\R_{n-1}$ with probability $q$. This proves Equations \eqref{eq:Orec}, \eqref{eq:Lrec}, \eqref{eq:Rrec}.

The table of type $\I_n$ is also straightforward. For table $\I_n$ we place the diner in seat $n-2$, which results in tables of type $\I_{n-3}$ and $\Lf_2$ with probability $p$, or in tables of type $\R_{n-3}$ and $\I_2$ with probability $q$. Thus,
\begin{align*}
 \E_B(\I_n) &= p(\E_B(\I_{n-3}) + \E_B(\Lf_2)) + q(\E_B(\R_{n-3}) + \E_B(\I_2)),\\
  &= p\E_B(\I_{n-3}) + q\E_B(\R_{n-3}) + p\E_B(\Lf_2) + q\E_B(\I_2), \\
  &= \E_B(\R_{n-2}) + \E_B(\Lf_3),
\end{align*}
where the final equality comes from Equations \eqref{eq:Lrec} and \eqref{eq:Rrec}. Thus we have established Equation \eqref{eq:Irec} and the proposition follows.
\end{proof}

Given the recurrences in Proposition \ref{prp:Erecs}, a straightforward induction argument (that we leave to the reader) proves the following intuitive properties for our expectations.

\begin{lemma}[Intuitive inequalities]\label{lem:intuitive}
Expectations increase monotonically, i.e., $0\leq \E_B(\mathcal{T}_n) \leq \E_B(\mathcal{T}_{n+1})$ for all $n\geq 0$. Furthermore, for each $n\geq 0$,
\begin{equation}\label{ineq:chain}
 \E_B(\Ou_n) \leq \E_B(\R_n) \leq \E_B(\Lf_n) \leq \E_B(\I_n).
\end{equation}
\end{lemma}

Now, for each table type $\mathcal{T}_n$, let $T(z)$ (in Roman font) denote the generating function for the sequence of expectations $\E_B(\mathcal{T}_n)$, i.e.,
\[
 T(z) = \sum_{n\geq 1} \E_B(\mathcal{T}_n) z^n = \E_B(\mathcal{T}_1) z + \E_B(\mathcal{T}_2) z^2+ \E_B(\mathcal{T}_3) z^3 + \cdots .
\]

We can use Proposition \ref{prp:Erecs} to express these generating functions in terms of each other with standard ``generatingfunctionology'' techniques \cite{Wilf}. We simply multiply both sides of each expression by $z^n$ and sum over all $n$ for which the statements are true. In this way, the recurrences become the following identities for the correspoding generating functions:
\begin{align*}
 C(z) &= z( pL(z) + qR(z) ),\\
 O(z) &= z( pO(z) + qR(z) ) = z( pL(z) + qO(z) ),\\
 L(z) &= z( pL(z) + qI(z) ),\\
 R(z) &= z( pI(z) + qR(z) ), \mbox{ and }\\
 I(z) &= z(1+z+z^2) + z^2 R(z) + \frac{q(1+p)z^4}{1-z},
\end{align*}
where the final term for $I(z)$ follows because $\E_B(\Lf_3) = q(1+p)$.

In solving this system of equations it is convenient to define the rational function $G(z)$ as follows:
\begin{equation}\label{eq:G}
 G(z) = \frac{z(1-p^2z^3)}{(1-z)^2(1-pz)(1+pz+pz^2)}.
\end{equation}

\begin{corollary}[Expectation generating functions]\label{cor:Egf}
The expectation generating functions for the various table types are:
\begin{align*}
 C(z) &= pqz^2(2-z)G(z), \\
 O(z) &= pqz^2G(z), \\
 L(z) &= qz(1-qz)G(z), \\
 R(z) &= pz(1-pz)G(z), \\
 I(z) &= (1-pz)(1-qz)G(z).
\end{align*}
\end{corollary}

From these generating functions, we derive the formula $O(z) = z(pL(z)+qO(z))$, giving us the following recurrence which will be useful later:
 \begin{equation}
     \E_B(\Ou_n) = p\E_B(\Lf_{n-1}) + q\E_B(\Ou_{n-1}). \label{eq:Orec2}
 \end{equation}

The form of these generating functions allows for quick evaluation of the limiting expectations. The concept is that each expectation grows at most linearly with the size of the table, and we can use a method similar to finite differences to find its growth rate. To see how this works, suppose $f(z) = \sum_{k\geq 0} a_k z^k$ and further suppose $\lim_{z \to 1} (1-z)^2f(z)=L$ exists. We write $\Delta a_k = a_k-a_{k-1}$, so that $(1-z)f(z) = a_0 + \sum_{k\geq 1}\Delta a_kz^k$. From this, we get that $(1-z)^2f(z) = a_0(1-z) + \Delta a_1z + \sum_{k\geq 2}(\Delta a_k-\Delta a_{k-1})z^k$. The limit of this function as $z \rightarrow 1$ is $L=\Delta a_1 + \sum_{k\geq 2}(\Delta a_k-\Delta a_{k-1})$. But notice the $n$th partial sum of this series telescopes as follows:
\begin{align*}
\Delta a_1 + \sum_{k=2}^n (\Delta a_k-\Delta a_{k-1}) &= \Delta a_1 + (\Delta a_2-\Delta a_1) + (\Delta a_3-\Delta a_2) + \cdots \\
&\quad + (\Delta a_n-\Delta a_{n-1}) = \Delta a_n.
\end{align*}
Therefore, $L = \lim_{n\to \infty} \Delta a_n$. A standard epsilon-delta argument now shows that $a_n = n\cdot L + r_n$, where $r_n/n \to 0$, i.e., $\lim_{n\to \infty} a_n/n=L$ as well.

Applying this observation to our generating function $G(z)$, we have $L=\lim_{z \to 1} (1-z)^2G(z) = \frac{1+p}{1+2p}$. The expressions in Corollary \ref{cor:Egf} are all multiples of $G(z)$ that go to $pq$ as $z \rightarrow 1$. This implies the following result, which says for large $n$, the expected proportion of napkinless diners is the same for each table type, and provides a formula for this expected proportion. We also remark that once we prove Algorithm B is optimal, this will establish the enumerative claim in Theorem \ref{thm:main}.

\begin{corollary}[Limiting expectations]\label{cor:limitE}
For any table type $\mathcal{T}_n$, we have following limiting proportion of napkinless diners:
\[
 \lim_{n\to \infty} \frac{\E_B(\mathcal{T}_n)}{n} =  pq\frac{1+p}{1+2p}.
\]
\end{corollary}

This corollary can also be understood intuitively as follows. Focus on the seats the maitre d' uses to set a wide trap. A trap has two empty seats with a left-reaching diner: $--\tl$. However, there might be a number of right-reaching diners in a row before a diner first reaches left: $--\tl\tr\tr\cdots\tr$. The probability of at least $k$ letters $\tr$ in this run of seats is $q^k$. Thus, on a large table the expected number of seats per trap is approximately
\[
 2 + (1+q+q^2+q^3+\cdots) = 2 + \frac{1}{1-q} = 2+\frac{1}{p}.
\]
Here's another way of understanding this quantity: a probability $p$ of reaching left can be interpreted as a 1 in $\frac{1}{p}$ chance, implying that on average, it takes $\frac{1}{p}$ attempts to get a diner to reach left. After a diner reaches left, 2 seats are left empty in order to set a trap, so the expected number of seats per trap is $2+\frac{1}{p}$.

These traps fail only when two consecutive diners reach left, so the probability of success for a trap is $1-p^2$. The proportion of napkinless diners, a.k.a., the number of successful traps divided by the number of seats, is 
\[
\frac{ \mbox{(number of successful traps)/(number of traps) } }{ \mbox{ (number of seats)/(number of traps) } },
\]
which by our estimates above is approximately $\frac{1-p^2}{2+\frac{1}{p}}= pq\frac{1+p}{1+2p}$ on a large table.

\section{Proof of optimality}

We are now ready to prove Algorithm B is an optimal strategy. To assist in the proof of optimality, we prove two subtle, technical lemmas for expectations. These could be proved using Proposition \ref{prp:Erecs} and induction, but we give tidy proofs of these facts using our generating functions.

\begin{lemma}[Subtle identity]\label{lem:technical}
 For all $n\geq 5$,
 \begin{equation}\label{eq:Lsubtle}
  p\E_B(\Lf_3) + p\E_B(\Lf_{n-3}) + q\E_B(\Lf_{n-1}) = \E_B(\Lf_n) + p^{n-3}q(1-2p).
 \end{equation}
\end{lemma}

\begin{proof}
To prove Equation \eqref{eq:Lsubtle}, let $H(z)$ denote the generating function for the quantity $p\E_B(\Lf_3) + p\E_B(\Lf_{n-3}) + q\E_B(\Lf_{n-1})-\E_B(\Lf_n) - p^{n-3}q(1-2p) = h_n$, defined for $n \geq 5$. We have $\E_B(\Lf_2) = q$, $\E_B(\Lf_3) = q(1+p)$, and $\E_B(\Lf_4) = q(1+p+p^2)$, so
\begin{align*}
 H(z) &= \sum_{n\geq 5} (\mbox{\small{ $p\E_B(\Lf_3)+p\E_B(\Lf_{n-3})+ q\E_B(\Lf_{n-1})-\E_B(\Lf_n)- p^{n-3}q(1-2p)$ }} )z^n,\\
  &= \frac{ pq(1+p)z^5 }{(1-z)}  + pz^3L(z) + qz( L(z) - qz^2-q(1+p)z^3 ) \\
  & \quad - (L(z) - qz^2- q(1+p)z^3 - q(1+p+p^2)z^4) - \frac{p^2q(1-2p)z^5}{(1-pz)}.
\end{align*}
With careful algebra (or the help of a computer algebra system) we see this expression simplifies to zero, so $h_n=0$ for all $n\geq 5$, which completes the proof.
\end{proof}

\begin{lemma}[Subtle inequality]\label{lem:technical2}
For all $n\geq 3$, we have $p^2q \leq \E_B(\R_n) - \E_B(\R_{n-1}) \leq pq$.
\end{lemma}

\begin{proof}
To prove our inequality, define $K(z)$ to be the generating function for the sequence of first differences, $k_n = \E_B(\R_n) - \E_B(\R_{n-1})$. We find 
\[
 K(z) = (1-z)R(z)= \frac{pz^2-p^3z^5}{1-qz-pz^3},
\]
which has a denominator of $1-qz-pz^3$ and a numerator of degree $5$, implying that for $n > 5$, $k_n = qk_{n-1} + pk_{n-3}$, which is a weighted average of $k_{n-1}$ and $k_{n-3}$ with positive weights. This implies that $k_n$ is between $p^2q$ and $pq$ if both $k_{n-1}$ and $k_{n-3}$ are between $p^2q$ and $pq$. We have $p^2q \leq k_i \leq pq$ for $i=3,4,5$ as base cases, so the lemma follows by induction. 
\end{proof}

To prove Algorithm B is an optimal strategy, we need to show that $\E_B(\mathcal{T}_n) = \E(\mathcal{T}_n)$ for all table types $\mathcal{T}_n$. Given the expressions for $\E(\mathcal{T}_n)$ in Equations \eqref{eq:ELmax}--\eqref{eq:EImax}, this boils down to establishing the inequalities in the following proposition. The case of circular tables then follows because Equations \eqref{eq:ECmax} and \eqref{eq:Crec} give the same recurrence for $\E(\C_n)$ as for $\E_B(\C_n)$.

\begin{proposition}[Local optimality]
The following inequalities hold for $1\leq i\leq n$:
\begin{align*}
\E_B(\Lf_n) &\geq \E_B(\Ou_i) + \E_B(\Lf_{n-i+1}),\\
\E_B(\R_n) &\geq \E_B(\R_i)+\E_B(\Ou_{n-i+1}), \\
\E_B(\Ou_n) &\geq \E_B(\Ou_i) + \E_B(\Ou_{n-i+1}), \mbox{ and } \\
\E_B(\I_n) &\geq \E_B(\R_i) + \E_B(\Lf_{n-i+1}).
\end{align*}
\end{proposition}

\begin{proof}
The proof is by induction on $n$. The base cases are true as we found in discovering Algorithm B.

Now for induction suppose that $\E(\mathcal{T}_k) = \E_B(\mathcal{T}_k)$ for all table types and all $k \leq n-1$.

For $\Lf_n$, we find, for any $2\leq i \leq n-1$:
\begin{align*}
 \E_B(\Lf_n) &= p\E_B(\Lf_{n-1}) + q\E_B(\I_{n-1}),\\
  &\geq p\E_B(\Ou_{i-1}) + p\E_B(\Lf_{n-i+1})\\
  & \quad + q\E_B(\R_{i-1})+q\E_B(\Lf_{n-i+1}),\\
  &= \E_B(\Ou_i) + \E_B(\Lf_{n-i+1}).
\end{align*}
The first and last equalities follow from Proposition \ref{prp:Erecs} and the inequality is by induction for tables $\Lf_{n-1}$ and $\I_{n-1}$. If $i=1$, $\E_B(\Ou_1) = 0$, so the claim is the trivial statement that $\E_B(\Lf_n) \geq \E_B(\Lf_n)$. If $i = n$, $\E_B(\Lf_1)=0$, so we have $\E_B(\Ou_n) + \E_B(\Lf_1) = \E_B(\Ou_n) \leq \E_B(\Lf_n)$ by Lemma \eqref{lem:intuitive}.

The argument for $\R_n$ is the same upon swapping $\Lf$ for $\R$, $i$ for $n-i+1$, and $p$ for $q$, and the case of $\Ou_n$ is also very similar.

The most interesting case is for the inner-facing table. For $1\leq i \leq n-4$,
\begin{align*}
 \E_B(\I_n) &= p\E_B(\I_{n-3}) + q\E_B(\R_{n-3}) + \E_B(\Lf_3), \\
  &= p\E_B(\I_{n-3}) + q(\E_B(\R_{n-3}) + \E_B(\Lf_3)) + p\E(\Lf_3) \\
  &= p\E_B(\I_{n-3}) + q\E_B(\I_{n-1}) + p\E_B(\Lf_3),\\
  &\geq p\E_B(\R_i) + p\E_B(\Lf_{n-i-2}) \\
  &\quad + q\E_B(\R_i) + q\E_B(\Lf_{n-i}) + p\E(\Lf_3),\\
  &=\E_B(\R_i) + p\E_B(\Lf_{n-i-2}) + q\E_B(\Lf_{n-i}) + p\E_B(\Lf_3), \\
  &=\E_B(\R_i) + \E_B(\Lf_{n-i+1}) + p^{n-i-2}q(1-2p), \\
  &\geq \E_B(\R_i) + \E_B(\Lf_{n-i+1}).
\end{align*}
Here we use the recurrences in Proposition \ref{prp:Erecs}, induction, as well as the ``subtle'' identity of Equation~\eqref{eq:Lsubtle} in the final equation. The case $i=n-3$ is $\E_B(\R_{n-3})+\E_B(\Lf_4) \leq \E_B(\R_{n-2})+\E_B(\Lf_4)-p^2q = \E_B(\R_{n-2})+\E_B(\Lf_3) = \E_B(\I_n)$ by Lemma \ref{lem:technical2}. For $i=n-2$, we have $\E_B(\R_{n-2}) + \E_B(\Lf_3) = \E_B(\I_n)$. With $i=n-1$, $\E_B(\R_{n-1})+\E_B(\Lf_2) \leq \E_B(\R_{n-2})+\E_B(\Lf_2)+pq = \E_B(\R_{n-2})+\E_B(\Lf_3) = \E_B(\I_n)$, again by Lemma \ref{lem:technical2}. Finally, for $i=n$ we have $\E_B(\R_n)+\E_B(\Lf_1) = \E_B(\R_n)$. This is at most $\E_B(\I_n)$ by Lemma \ref{lem:intuitive}. Thus, $\mathbb{E}_\mathrm{B}(\I_n) \geq \mathbb{E}_\mathrm{B}(\R_i) + \mathbb{E}_\mathrm{B}(\Lf_{n-i+1})$ for $1 \leq i \leq n$.
\end{proof}

Using equation \eqref{eq:Orec}--\eqref{eq:Rrec} and \eqref{eq:Orec2}, we can use the inequalities derived above to get:

\[\E_B(\Lf_n) \geq p(\E_B(\Lf_{i-1}) + \E_B(\Lf_{n-i})) + q(\E_B(\Ou_{i-1}) + \E_B(\I_{n-i})) \; :\; 1 \leq i < n \]
\[\E_B(\R_n) \geq p(\E_B(\I_{i-1}) + \E_B(\Ou_{n-i})) + q(\E_B(\R_{i-1}) + \E_B(\R_{n-i})) \; :\; 1 < i \leq n\]
\[\E_B(\Ou_k) \geq p(\E_B(\Lf_{i-1}) + \E_B(\Ou_{n-i})) + q(\E_B(\Ou_{i-1}) + \E_B(\R_{n-i})) \; :\; 1 \leq i \leq n\]
\[\E_B(\I_n) \geq  p(\E_B(\I_{i-1}) + \E_B(\Lf_{n-i})) + q(\E_B(\R_{i-1}) + \E_B(\I_{n-i})) \; :\; 1 < i < n\]

The form of these inequalities resembles Equations \eqref{eq:ELmax}--\eqref{eq:EImax}. Using those equations, as well as \eqref{eq:ECmax}, \eqref{eq:Crec}, Lemma \ref{lem:intuitive}, and the fact that $\E_B(\mathcal{T}_n) \leq \E(\mathcal{T}_n)$, one can show inductively that $\E(\mathcal{T}_n) = \E_B(\mathcal{T}_n)$ for all $n$. This proves the optimality of Algorithm B claimed in Theorem \ref{thm:main}.

\section{Probability generating function}

We have accomplished the task of proving Theorem 1, i.e., we have shown Algorithm B maximizes the expected number of napkinless diners, and we have provided a formula for the limiting expectation in terms of $p$. However, we can do even better and describe the entire probability distribution. That is, for table type $\mathcal{T}_n$, write
\[
 T_n(x) = \sum_{\sigma \in \{\tl,\tr\}^n} \pr(\sigma) x^{\nu(\sigma)} = \sum_{k\geq 0} \pr(n,k)x^k,
\]
where $\pr(n,k)$ denotes the probability that Algorithm B results in $k$ napkinless diners on a table of type $\mathcal{T}_n$.
We bundle up each of these polynomials into a bivariate generating function that has the full distribution of napkinless probabilities inside of it:
\[
 T(x,z) = \sum_{n\geq 1} T_n(x) z^n.
\]

The same idea behind the proof of Proposition \ref{prp:Erecs} can be used to find a system of recurrences for the polynomials $T_n(x)$, which in turn gives a solvable linear system for the $T(x,z)$. In particular, here is the expression we get for $C(x,z)$:
\begin{align}
 C(x,z) &= \frac{ z( 1 + pq(2x-1)z^2 - p^2z^3 - p^3q(x-1)z^5 )}{(1-pz)(1-qz - p(p+qx)z^3 -p^2q(x-1)z^4) } \label{eq:Cgf} \\
 &= z + z^2 + (p^2 +q^2 + 2pqx)z^3 + \cdots. \nonumber
\end{align}
As outlined in \cite[III.2.1]{FS}, not only the expectation but \emph{all} of the higher moments for the distribution are available to us once we have an explicit formula for $C(x,z)$. Using techniques of that book or \cite[Chapter 5]{Wilf}, we find 
\[
 \lim_{n\to \infty} \frac{\V(\C_n)}{n} =  \frac{pq(1+p)(1-p-2p^2+3p^3+8p^4)}{(1+2p)^3},
\]
where $\V(\C_n)$ denotes the variance of the distribution of $\nu(\sigma)$ on $\C_n$.

\begin{figure}
\includegraphics[width=15cm]{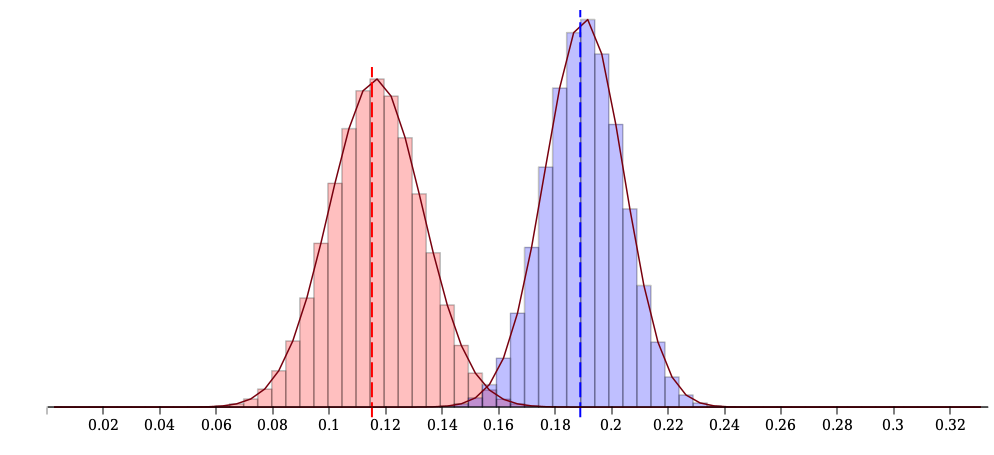}
\caption{Histograms representing the distribution of $\nu(\sigma)$ on a circular table with $n=200$ seats with $p=0.1543$ (in red) and $p=0.4554$ (in blue).}\label{fig:histograms}  
\end{figure}

By extracting coefficients $C_n(x)$ from this generating function, we can plot histograms corresponding to the distribution of napkinless diners. In Figure \ref{fig:histograms} we see the distribution of napkinless diners for $n=200$ seats (that's a big table!), and two different choices of probability: $p=0.1543$ (in red, with an expected $11.51\%$ napkinless) and $p=0.4554$ (in blue, with an expected $18.89\%$ napkinless). The higher value of $p$ maximizes the limiting proportion of napkinless diners, while the lower value of $p$ maximizes the limiting (normalized) variance of the number of napkinless diners.

\section{Final remarks and variants}

We are happy to have added another chapter to the saga of the malicious maitre d', but the story does not have to end here. We finish by discussing two potential directions for further study.

\subsection{French diners}

 There are several avenues for further investigation and generalization mentioned in \cite{APSTenner} and \cite{APSToal}; one of these is the problem of ``French diners'' as considered in \cite{Eriksen}. A French diner is one who will only take a napkin if it is their preferred napkin.\footnote{As mentioned in \cite{CP}, this modification was suggested by French mathematician Sylvie Corteel, who said if the ``proper'' napkin was not available, a French diner would simply cross their arms and refuse to eat.}  Let $f$ denote the proportion of French diners and assume that Frenchness and napkin preference are uncorrelated. For Conway's napkin problem of totally random seating, Eriksen found (in the $p=1/2$ case) the limiting proportion of napkinless diners to be 
\[
 \left( 1-\frac{ e^{\frac{1-f}{2}} -1}{1-f} \right)^2.
\]
This can be generalized using Eriksen's techniques to
\[
\left(\frac{q}{p}-\frac{e^{(1-f)q}-1}{1-f}\right)\left(\frac{p}{q}-\frac{e^{(1-f)p}-1}{1-f}\right)
\]
for general $p$. In the limiting case of $f = 1$, the expected proportion of napkinless diners is $pq$.

How does the proportion of French diners affect the optimal strategy for the adaptive maitre d'? We have started investigating this computationally using the approaches outlined in this paper. Here are our observations so far (we assume $0 < p \leq \frac{1}{2}$ throughout):
\begin{itemize}
\item Napkin shunning (by following step B1) chooses an optimal next seat when it is possible to do.
\item For $f = 1$, every strategy is optimal, and the expected proportion of napkinless diners on $\C_n$ for $n \geq 2$ is exactly $pq$. (We can prove this by noting that each napkin has a probability $pq$ of not being taken regardless of the strategy.) 
\item For $f < \frac{p}{q(1+p)}$, Algorithm B chooses the next seat optimally in large intervals. The limiting expected proportion of napkinless diners is $pq\frac{1+p+p^2f+pqf^2}{1+2p}$.
\item For $\frac{p}{q(1+p)} < f < 1$, one strategy that chooses an optimal next seat in large intervals is ``napkin shunning with small traps,'' i.e., it is the same as Algorithm B except we set one-seat-wide traps on the right. The limiting expected proportion of napkinless diners is $pq\frac{1+pf}{1+p}$.
\item For $f = \frac{p}{q(1+p)}$, the trap widths in the optimal strategy for large inner-facing tables follow a period-3 pattern: one-seat-wide traps in inner-facing tables of lengths 0 or 1 mod 3 and two-seat-wide traps in inner-facing tables of length congruent to 2 mod 3. The limiting expected proportion of napkinless diners is $\frac{p}{(1+p)^2}$. 
\item In the special case of $p = \frac{1}{2}$ and $f = \frac{2}{3}$, seating a diner in any non-endpoint seat in an inner-facing table with at least 3 seats is optimal. The expected proportion of napkinless diners on $\C_n$ for $n \geq 3$ is exactly $\frac{2}{9}$, with no lower-order terms!
\end{itemize}

\subsection{Sparse tables}

Another line of future investigation, which we believe is very intriguing, is to consider what happens when the number of diners is less than the number of seats. The approach of \cite{Eriksen} can be used to investigate Conway's version of this problem; if we let $d = \frac{n}{m}$ and let both $m$ and $n$ approach infinity while $d$ stays constant, we have that the limiting proportion of napkinless diners is \[\frac{(1-pe^{dq})(1-qe^{dp})-(1-d)pq}{dpq}.\] We close with some remarks about the adaptive ma\^itre d' version of this problem.

Suppose $m>n$ is the number of seats at the table. We see Algorithm B is clearly suboptimal since napkinless diners are among the last to sit. While Algorithm B may have created several intervals with fewer napkins than empty seats, the maitre d' might run out of diners before forcing anybody to become napkinless. At the same time, the maximum expected number of napkinless diners is always at least $\E_B(\C_{n})$ (so extra space cannot be a disadvantage). This is because we can leave a gap of $n-2$ seats between the first two diners (on the side that the first diner chooses) and use Algorithm B to fill in these seats, ignoring the other $m-n$ seats between diners 1 and 2. We can show that this results in an expectation of $\E_B(\C_{n})$ napkinless diners.

Here is our description of a promising strategy for $m = \infty$, i.e., an infinite straight table. This strategy is similar to Winkler's trap-setting strategy for $m = n$ in \cite[p. 25]{Winkler} and can be applied on sufficiently large finite tables. After seating the first diner, the strategy loops through the following three steps until the maitre d' runs out of diners:
\begin{itemize}
    \item Set one-seat-wide ``Winkler-style" traps until a trap succeeds.
    \item Seat a diner in the trapped seat.
    \item Reset by seating a diner at least $2n$ seats away from any other seated diner.
\end{itemize}
We can show that the limiting expected proportion of napkinless diners is $pq$ for this strategy and, with a little more work, that this limiting proportion is optimal. This strategy is optimal for a fixed value of $n\leq 7$ diners, and we previously believed that this strategy is optimal for all $n$. However, for $n = 8$, one chooses the unique optimal next seat by setting a 7-seat-wide trap. Clearly further investigation is required.

The contrast between Algorithm B and the strategy given above illustrates the tradeoff between saving space and saving diners. Setting a one-seat-wide trap is diner-efficient: if the trap succeeds, we only have to seat one diner in that trap to exploit it. On the other hand, setting a two-seat-wide trap is space-efficient: we need to seat two diners in a trap to exploit it (as opposed to one diner for a one-seat-wide trap), but these traps fail less frequently and therefore waste less space. We believe that this tradeoff will be a major part of the dynamics in this new variant.

\bibliographystyle{plain}
\bibliography{bibliography2}

\vfill\eject

\end{document}